# A Proof, Based on the Euler Sum Acceleration, of the Recovery of an Exponential (Geometric) Rate of Convergence for the Fourier Series of a Function with Gibbs Phenomenon

John P. Boyd[*]
Department of Atmospheric, Oceanic
and Space Science and Laboratory for Scientific Computation,
University of Michigan, 2455 Hayward Avenue, Ann Arbor MI 48109
jpboyd@engin.umich.edu;
http://www.engin.umich.edu:/~ jpboyd/

March 26, 2010

## Abstract

When a function $f(x)$ is singular at a point $x_s$ on the real axis, its Fourier series, when truncated at the $N$-th term, gives a pointwise error of only $O(1/N)$ over the entire real axis. Such singularities spontaneously arise as "fronts" in meteorology and oceanography and "shocks" in other branches of fluid mechanics. It has been previously shown that it is possible to recover an exponential rate of convergence at all points away from the singularity in the sense that $|f(x) - f_N^\sigma(x)| \sim O(\exp(-q(x)N))$ where $f_N^\sigma(x)$ is the result of applying a filter or summability method to the partial sum $f_N(x)$ and $q(x)$ is a proportionality constant that is a function of $d(x) \equiv |x - x_s|$, the distance from $x$ to the singularity. Here we give an elementary proof of great generality using conformal mapping in a dummy variable $z$; this is equivalent to applying the Euler acceleration. We show that $q(x) \approx \log(\cos(d(x)/2))$ for the Euler filter when the Fourier period is $2\pi$. More sophisticated filters can increase $q(x)$, but the Euler filter is simplest. We can also correct recently published claims that only a root-exponential rate of convergence can be recovered for filters of compact support such as the Euler acceleration and the Erfc-Log filter.

## 1 Introduction

Suppose that $f(x)$ is a function with a Fourier series which converges poorly because of a singularity at single point $x_x \in [-\pi, \pi]$:

$$f(x) = \sum_{n=-\infty}^{\infty} c_n \exp(inx) \qquad (1)$$

---
[*]

(We shall assume that the coefficients $c_n$ are all bounded by a constant $C$.) "Converges poorly" means either a "subgeometric" or "algebraic" rate of convergence in the language of [8]: it must be impossible to find a bound of the form

$$|c_n| \leq \text{constant} \exp(-q|n|) \qquad (2)$$

for any finite real $q > 0$; if such a bound is possible, then Theorem 2 below shows that $f(x)$ must be singularity-free everywhere within the strip $\Im(x) < q$, which contradicts the assumption that there is a singularity on the real axis.

For example, if $f(x)$ has a jump singularity at $x = x_s$, then the $c_n$ decay only a $O(1/|n|)$. The $N$-th partial sum

$$f_N(x) = \sum_{n=-N}^{N} c_n \exp(inx) \qquad (3)$$

yields an approximating only $f(x)$ of only $O(1/N)$ over the entire real axis — worse in the Gibbs boundary layers of thickness $O(1/N)$ around the discontinuity (and its copies at $x = x_s + 2m\pi$ for all integer $m$) where the error is $O(1)$ [38].

It has been known for several decades that the rate of convergence of convergence of a slowly converging Fourier series can be greatly improved by applying a linear filter.

**Definition 1 (Linear Sum Acceleration Filter)** *A linear sum acceleration filter approximates $f(x)$, the limit of the partial sums as $N \to \infty$, by the filtered partial sum*

$$f_N^\sigma \equiv \sum_{n=-N}^{N} \sigma(j/N) c_n \exp(inx) \qquad (4)$$

*where $\sigma(\theta)$ is symmetric with respect to $\theta = 0$, that is, $\sigma(-\theta) = \sigma(\theta)$ for all $\theta$. In plain language, a filter multiplies each coefficient of the spectral series by a numerical weight factor, $\sigma(j/N)$.*

**Definition 2 (Filter Order)** *A filter is said to be order $p$ if, for any function $f(x)$ which is periodic and analytic for all real $x$ (and therefore does not need acceleration),*

$$|f(x) - f_N^\sigma(x)| \leq \frac{\text{constant}}{N^p} \qquad (5)$$

Majda, McDonough and Osher [46] showed in 1978 that is possible to recover spectral or near-spectral accuracy by using an filter which is an exponential. However, their filter does not precisely satisfy the order conditions defined by later authors.

Eight years later, Boyd and Moore [13] showed that the Euler acceleration could recover *exponential* accuracy for the closely related problem of slowly-converging Hermite series.

In 1991, Vandeven [66] defined a filter to be of order $p$ as above. He gave necessary and sufficient conditions for a filter to be of order $p$. In his Theorem 3, he proved that for a jump discontinuity in $f(x)$

$$|f(x) - f_N^\sigma(x)| \leq \frac{\text{constant}}{d(x)^{p-1} N^{p-1}} ||u||_p \qquad (6)$$

where $d(x)$ is the distance from $x$ to the nearest singularity on the real $x$-axis. In his Theorem 4, he proved that recovery of an *exponential* rate of convergence is possible. However, the rate of convergence is *subgeometric* in the sense that the bound is proportional to $\exp(\text{constant} N^r$ for some $r < 1$ [8]. On a fixed subinterval, $|x| \geq 1$ with the singularity at the origin, he proved that his filter could recover an error no worse than $N^\beta \exp(N^{1/4} \{\log(C) - 1/2 \log(N)\})$ where $\beta > 0$ and $C$ are positive constants independent of $N$. (He also proved slower rates of convergence for subintervals that converge towards the singularity, but we shall not repeat them here.)

Four years later, Boyd applied the Euler sum acceleration [5] to the piecewise linear or "sawtooth" function. For conformity with the convention used here of a singularity at the origin, we shall replace his sawtooth function by the "shifted sawtooth" defined by $Sws(x) = Sw(x + \pi)$:

$$Sws(x) \equiv \begin{cases} x - \pi & x \in [0, 2\pi] \\ Sws(x + 2\pi m) & \text{otherwise}, m = \text{integer} \end{cases} \quad (7)$$

which has the very slowly converging Fourier series

$$Sws(x) \equiv -2 \sum_{n=1}^{\infty} \frac{1}{n} \sin(nx) \quad \forall x \quad (8)$$

This is the prototype for function with discontinuities, which aerodynamicists call "shocks" and meteorologists call "fronts". He proved that for this example, the error at fixed $x$ of the Euler-accelerated series is proportional to $\exp(-\log(\cos(x/2))N)$. property, it follows that In [6], he developed an improved acceleration dubbed the "Erfclog" filter which is justified by a combination of an asymptotic approximation to Vandeven's filter with the "lagged-Euler" improvement of his earlier article [5].

Tadmor and Tanner [58, 59] using a different line of argument proved that with their filter, which they dubbed simply the "adaptive" filter, it was possible to obain what they called "root-exponential" convergence, that is, an error proportional to the exponential of $\sqrt{N}$. This subgeometric rate of convergence is asymptotically inferior to the Euler acceleration, which is an exponential *linear* in $N$ [5, 6].

Tanner [60] introduced a new filter which he dubs the Hermite Distributed Approximating Function (HDAF) and proved that it recovers an accuracy proportional to $\exp(-q(x)N)$. In his numerical experiments, this was a little better than the Erfclog filter of Boyd [6]; both were greatly superior to the adaptive filter of Tadmor and Tanner [58, 59].

He remarks that "Other examples of noncompactly supported filters include are the adaptive Erfc and Erfc-log filters". This is a little confusing because the Erfclog filter function $\sigma^{Erfclog}(\theta)$ is identically zero for all $|\theta| > 1$ and therefore is of compact support in Fourier space. However, Tanner applies the convolution theorem to express filtering as a convolution integral in $x$ rather than a simple product of Fourier terms. The kernel of the convolution integral is the "mollifier", which is just the Fourier transform of the filter function $\sigma$. The mollifier (Fourier transform) of the Erfclog filter is not of compact support, so his assertion is not contradicted by the exponential convergence of the Erfclog filter [6] if "compact support" is restricted to "compact support of the mollifier".

Both Tanner's HDAF filter and the Erfclog filter of Boyd [6] contain an order parameter $p$, and both are best applied with *spatially varying order*. Although the proportionalities are different, for both the rate of convergence is best when the order $p$ is *large* in smooth areas and decreases to zero at the point where $f(x)$ is singular.

The Euler acceleration is non-adaptive in the sense that the same filter weights are applied for small $x$. Although the proportionality factor $q(x)$ in error bound proportional to $\exp(-q(x)N)$ is not as large as for the adaptive HDAF and Erfclog methods, it is nevertheless true that the Euler acceleration gives an exponential (geometric) rate of convergence. It also permits a theoretical analysis which requires only undergraduate mathematics.

Table 1: List of Symbols

| Symbol | Name | Comments |
|---|---|---|
| $S$ | target sum | |
| $S^{Abel}(z)$ | Abel extension | |
| Sws | shifted sawtooth function | Sws $= x - \pi, x \in [0, 2\pi]$ |
| $z$ | | $z = 1$ gives the original sum |
| $\zeta$ | transformed coordinate | $\zeta = 1$ gives the original sum |
| $\mu$ | asymptotic rate of convergence | $S - S_N \propto \exp(-N\mu)$ |
| $\rho$ | convergence factor | equals the radius of convergence of $\zeta$ power series; also $\rho = \exp(\mu)$ |

Recovery of spectral accuracy in the presence of shocks and other singularities has become a "big business". A partial list of other efforts to accelerate Fourier series include [31, 9, 32, 30, 29, 28, 27, 24, 25, 23, 22, 21, 20, 17, 36, 19, 40, 39, 41, 42, 43, 44, 51, 45, 2, 3, 47, 48, 61, 54, 52, 70, 71, 57, 67, 68, 11, 10] It therefore seems worthwhile to give a simple proof that this "spectral recovery" is founded on stone instead of sand.

## 2 Inflating a Summation to a Power Series in a Dummy Variable $z$ Followed by Conformal Mapping

A common procedure for accelerating a slowly convergent series is the following three-step procedure:

1. Inflate the series to a function of a dummy variable $z$ by multiplying the $n$-th term of the series $z^n$

2. Apply a conformal mapping by replacing by a new $\zeta$ where

$$z = \mathcal{Z}(\zeta) \qquad (9)$$

with the mapping chosen to be an analytic function such that, among other properties, (i) $\mathcal{Z}(1) = 1$ and (ii) $\mathcal{Z}(z) \propto \zeta$ for small $\zeta$.

3. Expand the $N$-th partial sum of the inflated series as a power series in $\zeta$; because of the requirement that $z \propto \zeta + O(\zeta^2)$ for small $\zeta$, the first $N$ terms in the $\zeta$ series are completely determined by only the first $N$ terms in the original series.

The accelerated approximation is then just the partial sum of the $\zeta$ series evaluated at $\zeta = z = 1$.

Reviews are given by [50, 55, 33]. Scraton [56] and Bruno and Reitich [16] give further improvements.

The inflated series is dubbed the "Abel extension" of the sum $S$ in [13]. Since Niels Henrik Abel died in 1829, the counter-intuitive idea of simplifying a problem by making it more complicated by converting the sum into a power series has a long history!

A concrete example is

$$S^{log2} \equiv \log(2) = \sum_{n=1}^{\infty} (-1)^{n+1}/n \qquad (10)$$

The inflated series or "Abel Extension" is

$$S^{log2,Abel}(z) = \log(1+z) = \sum_{n=1}^{\infty} \frac{(-1)^{n+1}}{n} z^n \qquad (11)$$

The Abel extension is logarithmically singular at $z = -1$. Consequently, its power series in $z$ has a unit radius of convergence. At $z = 1$, the point corresponding to the original series for $\log(2)$, the power series converges very slowly because this point is right on the boundary of the disk of convergence. And yet the inflated function $S^{log2,Abel}(z)$ is *not* singular at $z = 1$, the only value of $z$ that we really care about, but only at $z = -1$.

If the conformally-mapped function $S^{log2,Abel}(z(\zeta))$ is expanded as a power series in $\zeta$, the domain of convergence will again be a disk, but this time in the $\zeta$-plane. However, this $\zeta$-disk maps back to a region of *different shape* in the $z$-plane. If we choose the conformal map $m(\zeta)$ wisely so that $z = 1$ now lies *within* the region of convergence, then the power series in $\zeta$ will converge geometrically at $z = 1$ (to the original unflated sum) in the sense that error will decrease as $\exp(-qN)$ for some constant $q > 0$.

"Conformal" is an adjective that means, in the narrowest sense, "angle-preserving". We do not actually use the angle-preserving property here. Rather, we merely wish to choose the mapping function $m(z)$ to be an *analytic* function so that the power series in $\zeta$ has a large a radius of convergence as possible. This analyticity automatically implies that the mapping is conformal. It would perhaps be more logical to refer to this as the "analytic mapping acceleration", but the terminology "conformal mapping acceleration" has become standard.

One choice is to pick a simple rational transformation that maps the singular point, $\zeta = -1$, to $\zeta = \infty$.

$$z = \frac{\zeta}{2-\zeta} \qquad \leftrightarrow \qquad \zeta = \frac{2z}{1+z} \qquad (12)$$

A rational map which is linear in both the numerator and denominator polynomials is also known as a "Möbius transformation". It has the important property that all Möbius transformations map circles to circles. Although we shall not demonstrate it, it can be shown that this mapping generates the same accelerated series as applying the filter weights devised by Leonard Euler in the eighteenth century [13].

The circles of constant $\zeta$ are "equiconvergence contours" in the sense that the $n$-th term of the $\zeta$ power series will decay proportional to $(|\zeta|/\rho)^n$ where $\rho$ is the radius of convergence of the power series in the complex $\zeta$-plane. The following classical theorem, found in most calculus texts, provides a more precise description.

**Theorem 1 (Convergence of a power series)** *Suppose that the disk $|\zeta| < \rho$ is the largest disk centered on the origin in the complex $\zeta$-plane which is free of singularities of an analytic function $S(\zeta)$. Then the power series of $S(\zeta)$ converges everywhere within the disk and thus $\rho$ is the radius of convergence of the series. Let $S_N(\zeta)$ denote the partial sum of the power series up to and including the $N$-th term. Then*

$$|S(\zeta) - S_N(\zeta)| \leq constant \exp(-N \log(\rho - \epsilon)) \qquad (13)$$

*for all $N$ with a proportionality constant independent of $N$ where $\epsilon > 0$ is a constant that may be arbitrarily small.*

**Definition 3 (ASYMPTOTIC RATE OF GEOMETRIC CONVERGENCE )** *If a series has geometric convergence, that is, if an expansion has an exponential index of convergence $r = 1$ so that*

$$a_n \sim [\ ] \exp(-n\mu) \qquad (14)$$

Table 2: Summary of Euler Acceleration for the log(2) Series

This example is special in that it is possible to find a closed-form expression for the Abel Extension in $z$ and $z(\zeta)$ as given below; the method can be applied even if only a finite number of terms of the original series are known.

| $S^{log2}$ | $\log(2)$ | $\sum_{n=1}^{\infty} \frac{(-1)^{n+1}}{n}$ |
|---|---|---|
| $S^{log2,Abel}(z)$ | $\log(1+z)$ | $\sum_{n=1}^{\infty} \frac{(-1)^{n+1}}{n} z^n$ |
| $S^{log2,Abel}(z(\zeta))$ | $-\log(1-\zeta/2)$ | $\sum_{n=1}^{\infty} \frac{(-1)^{n+1}}{n} \left(\frac{\zeta}{2-\zeta}\right)^n$ $-\sum_{n=1}^{\infty} \frac{2}{n} \left(\frac{1}{2}\right)^n \zeta^n$ |

*where $a_n$ are the spectral coefficients, $\mu$ is a constant, and [ ] denotes unspecified factors that vary more slowly with $n$ than the exponential (such as $n^k$ for some k), then the ASYMPTOTIC RATE OF GEOMETRIC CONVERGENCE is $\mu$. Equivalently,*

$$\mu = \lim_{n\to\infty} \{-\log|a_n|/n\} \tag{15}$$

*This definition is meaningful only for geometrically converging series; it does not apply when the algebraic index of convergence is $< \infty$ nor when the exponential index of convergence $r < 1$.*

The function $(z[\zeta])$ will be singular at $\zeta = 2$, which is the image of $z = \infty$, if $\tilde{S}(z)$ is singular at $z = \infty$. This is usually the case. It follows that the usual situation is that the $\zeta$ series converges only for $|\zeta| \leq 2$; as discussed in the next section, the radius of convergence may be smaller if the Abel extension has singularities off the real $z$-axis. However, $\log(1+z)$ has singularities only at $z = -1$ and $z = \infty$; the transformed function is singular only at $\zeta = 2$ and $\zeta = \infty$:

$$S^{log2,Abel}(z(\zeta)) = -\log(1-\zeta/2) = -\sum_{n=1}^{\infty} \frac{2}{n} \left(\frac{1}{2}\right)^n \zeta^n \tag{16}$$

The convergence region is

$$|\zeta| \leq 2 \quad \leftrightarrow \quad \Re(z) \geq -1/2 \tag{17}$$

That is, the $\zeta$ series converges everywhere in the half-plane to the right of the line $\Re(z) = -1/2$. At $z = \zeta = 1$, the coefficients of the Eulerized series decrease as $1/2^n$, multiplied by some slower-than-exponential function of $n$. Equivalently, the $n$-th term decays proportionally to $\exp(-\log(2)|n|)$.

Although we shall not derive this well-known result here, the result of evaluating the $\zeta$ series at $\zeta = 1$ is equivalent to replacing the $N$-th partial sum of the original series $S_N = \sum_{n=0}^{\infty} a_n$ by

$$S_N^{Euler} = \sum_{n=0}^{\infty} a_n \sigma_E(n/(N+1)) \tag{18}$$

$$\sigma_E(0) = 1 \tag{19}$$

$$\sigma_E\left(\frac{j}{M+1}\right) = \sum_{k=j}^{M} \mu_{Mk}, \quad j = 1, 2, \ldots M$$

$$\sigma_E(1) = 0$$

where the "partial sum weights" are

$$\mu_{Mk} = \frac{M!}{2^M} \frac{1}{k!(M-k)!} \tag{20}$$

This strategy is very general and can be applied to series that contain a parameter as we shall do in the next section. Pearce [55] reviews many different conformal mappings that have been used in practice. Steven Weinberg, who won the Nobel Prize in physics for his contributions to the unification of the weak interactions with electrodynamics, used the mapping

$$z = \frac{4\zeta r(g)}{g(r(g) - \zeta)^2} \tag{21}$$

where $g$ is a parameter and $r(g)$ is such that $z(1) = 1$; this mapping was a little better than the Euler map (which he calls a "Möbius transformation") for his problem for which it is known that all the singularities of the Abel extension lie on the negative real $z$-axis [69].

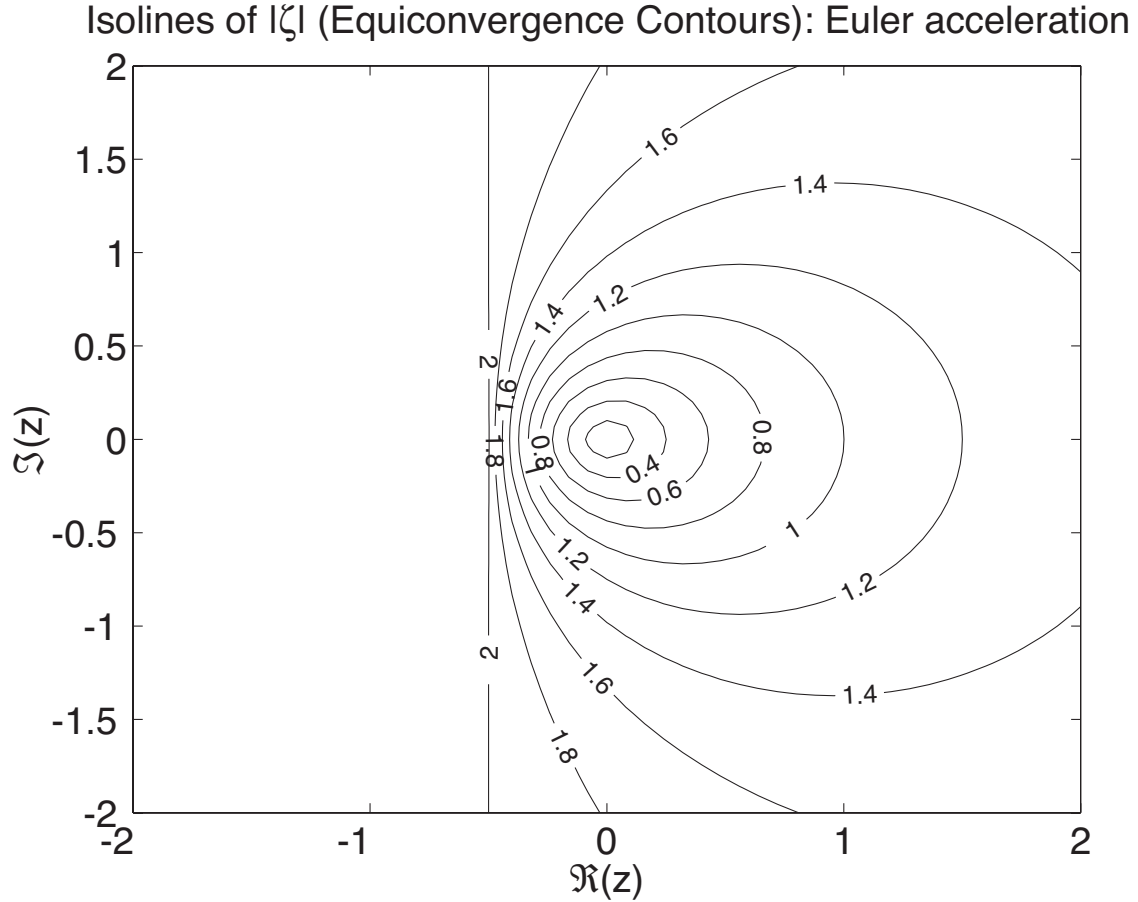

Figure 1: The Euler conformal map, a special case of the Möbius transformation, maps circles to circles. Each of the circular contours is the image of the circle of constant $|\zeta|$ in the $\zeta$-plane whose magnitude labels that contour. The radius of convergence $\rho$ of the power series in $\zeta$ is the label of the largest of the contours shown which is singularity-free in the $z$-plane

Table 3: Selected Applications of Euler Acceleration

| Application | Discipline | Reference |
|---|---|---|
| Blasius Flow | Fluid Mechanics | Boyd [7] |
| Stokes' Drag, Pair of Spheres | Hydrodynamics | Frost & Harper [26] |
| Oseen Drag: Flow Past a Sphere | Hydrodynamics | Hunter & Lee [37] Van Dyke [62] |
| Flow Between Rotating Disks | Hydrodynamics | Hoffman [35] Van Dyke [63] |
| Skin Friction, Parabolic Cylinder | Aerodynamics | Van Dyke [64, 65] Frost & Harper [26] |
| Standoff Distance, Blunt Body Shock | Aerodynamics | Van Dyke [64, 65] |
| Sawtooth Shock | Aerodynamics | Boyd [5] |
| Yoshida Jet | Oceanography | Boyd & Moore [49, 13] |
| Sinc Pseudospectral Differentiation | Numerical Analysis | Boyd [4] |
| Hyperasymptotic Improvement of the Asymptotic Series for the Stieltjes Function | Numerical Analysis | Olver [53] pp. 536-544 |
| Oscillatory Integrand | Numerical Analyis | Abramowitz & Stegun [1] pg. 22 |
| Born scattering series | physics | Weinberg [69] |

For the log(2) series, the only singularities are at $z = -1$ and $z = \infty$. Instead of mapping the finite singularity to infinity and thereby mapping $z = \infty$ to $\zeta = 2$ in the process, a better compromise is to map *both* singularities to $|\zeta| = 3$ by

$$z = \frac{2\zeta}{3 - \zeta} \qquad \leftrightarrow \zeta = \frac{3z}{2 + z} \qquad (22)$$

The terms of the $\zeta$ power series now decay in proportional to $(1/3)^n$ instead of $(1/2)^n$.

Other variants of the Möbius transformation are useful in special situations. If a power series in $z$ is a function of $z^2$ with convergence limited by singularities on the imaginary $z$-axis, the mapping $z = \sqrt{\frac{\zeta}{2-\zeta}}$, equivalent to $\zeta = 2z^2/(1 + z^2)$, yields a power series in $\zeta$ that will converge for all real $z$. Similarly, Boyd exploited a $C_3$ rotational symmetry in the complex $z$-plane for the Blasius function through the map $z = (\zeta/(2 - \zeta))^{1/3}$ which successfully extended Blasius' power series of 1908 into an expansion converging on the entire semi-infinite physical domain [7].

Rather than attempt to exhaustively review all the possibilities, we refer the reader to Pearce's review [55] and to the table of applications, Table 3.

The crucial point, though, is that because $z$ is linearly proportional to $\zeta$ for small $|\zeta|$, the first $N$ terms of the Eulerized series are completely determined by the first $N$ terms of the original series in $z$. Thus, it is completely unnecessary to know the *general* term in the series; *numerical computation of N terms is sufficient.*

# 3 Accelerating a Fourier Series

The goal is to improve the rate of convergence of

$$f(x) = \sum_{n=-\infty}^{\infty} c_n \exp(inx) \qquad (23)$$

One important complication compared to the previous section is that the sum (43) depends on the parameter $x$. We will find that the rate of convergence of the accelerated series varies strongly with $x$.

**Definition 4 (Pointwise and Global Accelerators)** *A pointwise accelerator for a Fourier (or other spectral series is one that treats each point in $x$ independently of the behavior of the sum at other $x$. A global accelerator uses information from all $x$ to reconstruct the sum at each $x$.*

**Proposition 1 (Failure of Pointwise Recombination Accelerators)** *All pointwise acceleration schemes, whether linear or nonlinear, are bound to fail sufficiently close to a jump discontinuity.*

We shall not offer a formal proof, but only an instructive example. The piecewise-linear or "sawtooth" function defined by () above is discontinuous at $x = 0$. However, its Fourier series is a *sine* series. The value of Sw($\epsilon$) for $\epsilon \ll 1$ is approximately $-\pi$, but all the terms in the sine series are *individually zero*. There is no way to combine zero terms to manufacture something nonzero.

Thus, the best solution to extract an accurate approximation near the discontinuity is "Gibbs Reprojection", which is to expand the Fourier series in a different basis of polynomials or radial basis functions. Gibbs Reprojection is a *global* acceleration method in the sense that information from the entire $x$ interval is used to reconstruct the function $f(x)$ at a given $x$ [34].

Nevertheless, the Euler method is successful everywhere except in the immediate vicinity of the discontinuity. For conformity with our convection here, we shall shift the sawtooth function defined in [8, 5]For example, the Abel extension of the sawtooth function is [5]

$$S^{Sws,Abel}(x, z) \equiv \begin{cases} x - \pi + 2\arctan\left(\frac{1-z}{1+z}\cot(x/2)\right) & x \in [0, 2\pi] \\ S^{Sws,Abel}(x + 2\pi m) & \text{otherwise}, m = \text{integer} \end{cases} \qquad (24)$$

The inflated sawtooth function is singular in the $z$-plane for a given $x$ at

$$z_s = \pm \exp(\pm ix) \qquad (25)$$

(Note that an overall minus sign was omitted in (2.5) of [5].) Boyd showed that, applying the Moebius conformal map, the series through the $N$-th term has an error proportional to

$$\text{Sws}(x) - f_N^s(x) \sim O\left(\cos^N\left(\frac{d(x)}{2}\right)/N\right) \qquad (26)$$

where $d(x)$ is the distance from $x$ to the nearest singularity of the sawtooth function. (Here, $d(x) = |x|$.)

However, the sawtooth function is not completely satisfactory because (i) other types of singularities besides discontinuities may arise and (ii) most functions have additional singularities off the real axis. In the next section, we shall analyze this general case.

# 4 Functions with Singularities in the Complex Plane

(As noted earlier, we assume that the coefficients $c_n$ are bounded by a constant $C$.) There is no loss of generality in assuming that the spatial period is $2\pi$ since a general period $P$ in a variable $y$ may be reduced to (43) by the trivial dilation $x = (2\pi/P)y$.

The most significant restriction we shall impose is that the Fourier series is slowly converging because of a single singularity on $x \in [-\pi, \pi]$. (We shall discuss multiple singularities later.) We shall not impose any restrictions on the type of singularity.) We shall *not* exclude the possibility that $f(x)$ has singularities off the real axis at $x = x_k, k = 1, 2, \ldots, n_s$ where $n_s$ is any integer between zero and infinity and $|\Im(x_k)|$ is nonzero for all $k$.

Without loss of generality, we shall assume that the singularity on the real axis between $x = -\pi$ and $x = \pi$ is located at $x = 0$. A singularity at a point $w = x_0$ in a coordinate $w$ can be translated into our standard form by defining $x = w - x_s$.

To inflate $f(x)$ to a power series in $z$, it is helpful to write

$$f(x) = f_1(x) + f_2(x) \tag{27}$$

where

$$f_1(x) = \sum_{n=0}^{\infty} c_n \exp(inx) \tag{28}$$

$$f_2(x) = \sum_{n=1}^{\infty} c_{-n} \exp(-inx) \tag{29}$$

Then the Abel extension is

$$f^{Abel}(x) = f_1^i(x) + f_2^i(x) \tag{30}$$

$$f_1^i(x, z) = \sum_{n=0}^{\infty} c_n z^n \exp(inx) \tag{31}$$

$$f_2^i(x, z) = \sum_{n=1}^{\infty} c_{-n} z^n \exp(-inx) \tag{32}$$

As in the previous section, we will apply the Euler acceleration to map $z$ to a new coordinate $z$ by a Möbius transformation and sum the first $N$ terms of the power series in $\zeta$ for some user-chosen $N$ at $z = \zeta = 1$ to recover $f(x)$. To determine the $x$-dependent rate of convergence, we first need several theorems and definitions.

**Definition 5 (Gevrey Space)** *A function $f(x)$ is in the periodic Gevrey space $G_{q,m}(-\pi, \pi)$ if*

$$||f||^2_{G_{q,m}(-\pi,\pi)} \equiv \sum_{n=-\infty}^{\infty} exp(2q(1 + |n|)) \left(1 + |n|^{2m}\right) |c_n|^2 < \infty \tag{33}$$

*where the $c_n$ are the Fourier coefficients of $f(x)$.*

Our assumption that $f(x)$ is singular at $x = 0$ is equivalent to assuming that $f(x)$ does *not* belong to any Gevrey space. However, in the proofs below, we shall consider some functions that hypothetically do belong to Gevrey spaces.

**Theorem 2** *If the Fourier coefficients of the series*

$$f(x) = \sum_{n=-\infty}^{\infty} c_n \exp(inx) \qquad (34)$$

*satisfy a bound of the form*

$$|c_n| \leq \text{constant} \exp(-q|n|) n^{2m} \qquad (35)$$

*for some positive constant $q$ and finite positive integer $m$, or equivalently belong to the Gevrey space $G_{q,m}(-\pi, \pi)$ then the sum of the Fourier series is spatially periodic,*

$$f(x + 2\pi) = f(x) \qquad \forall x \qquad (36)$$

*and is* ANALYTIC *everywhere within the strip*

$$|\Im(x)| < q \qquad (37)$$

*Noted on pg. 271 of [18]*

"Geometric convergence" means that the terms of a series can be bounded by those of the geometric series $1/(1-w) = \sum_{n=0}^{\infty} w^n$ for some $w$ where in this instance $w = \exp(-\mu)$. For complex values of $x$, the Fourier series for $f_1(x)$ and $f_2(x)$ can exhibit geometric convergence. The previous theorem then easily implies the following.

**Theorem 3** *The function $f_1(x) \equiv \sum_{n=0}^{\infty} c_n \exp(inx)$ has no singularities in the upper half-plane $\Im(x) > 0$ while $f_2(x) \equiv \sum_{n=1}^{\infty} c_{-n} \exp(-inx)$ is free of singularities in the lower half-plane.*

Proof. Let $x = \sigma + i\tau$. Then

$$f_1(\sigma + i\tau) = \sum_{n=0}^{\infty} c_n \exp(-n\tau) \exp(in\sigma) \qquad (38)$$

When $\tau = \Im(x) > 0$, the series converges geometricallly. Theorem 2 then implies that $f_1(x)$ must be analytic for all $\sigma$. The same reasoning applied to $f_2$ in the lower half-plane where $\exp(-inx) = \exp(-in\sigma)\exp(n\tau)$. QED

**Theorem 4 (Location of Singularities in $z$)** *Suppose $f(x)$ has a singularity at*

$$x_j = \sigma_j + i\tau_j, \qquad j = 0, 1, 2, \ldots \qquad (39)$$

*Then $f^{Abel}(x, z)$ has singularities at, for real $x$,*

$$\begin{aligned}
\tau_j < 0, &\qquad |z| = \exp(|\tau_j|), &\qquad \theta_j(x) = \sigma_j - x \\
\tau_j > 0, &\qquad |z| = \exp(\tau_j), &\qquad \theta_j(x) = x - \sigma_j \\
\tau_0 = 0, &\qquad |z| = 1, &\qquad \theta(x) = \pm x
\end{aligned} \qquad (40)$$

Proof. With

$$x = \sigma + i\tau, \qquad z = r\exp(i\theta) \tag{41}$$

$$f_1^i(x,z) = \sum_{n=0}^{\infty} c_n \exp(n[\log(r) - \tau]) \exp(in[\theta + \sigma]) \tag{42}$$

$$f_2^i(x,z) = \sum_{n=1}^{\infty} c_{-n} \exp(n[\log(r) + \tau]) \exp(-in[\theta - \sigma]) \tag{43}$$

Thus, $f_1^i(x,z)$ does not depend on $r$ and $\tau$ separately, but only on the combination $\log(r) + \tau$. Likewise, $f_1^i(x,z)$ does not depend on $\theta$ and $\sigma$ separately, but only on the combination $\theta + \sigma$. Similar remarks apply, except for a sign difference in the combinations, to $f_2^i(x,z)$.

It follows that if $f_1(x) = f_1^i(x, z = 1)$ has a singularity at $x = \sigma + i\tau$ where, as already proved, $\tau_j \le 0$, then

$$\log(r_j(\tau)) = \tau - \tau_j \quad \& \quad \theta_j(\sigma) = \sigma_j - \sigma \tag{44}$$

In particular, when $x$ is real and therefore $x = \sigma$, these become

$$|z_j| = \exp(|\tau_j|), \qquad \theta_j(x) = x - \sigma_j \; \forall \, \text{real} x, \; \tau_j \le 0 \tag{45}$$

Applying the same reasoning to $f_2^i$ proves the rest of the theorem. QED

Note that by assumption, there is only a single real singularity, that at $x = 0$; we will assign index 0 to this singularity. Note further that this singularity is usually the sum of singularities in $f_1(x)$ and $f_2(x)$, and thus the image of lone singularity at $x = 0$ is usually a *pair* of singularities on the unit circle in $z$; the exceptions are when $x = 0, \pm\pi, \pm 2\pi, \ldots$ where the singularities coincide in location in $z$.

If $f(x)$ is real for real $x$, which is unnecessary for our theory but is usually true in applications, then singularities off the real axis occur in complex conjugate pairs. Each pair generates a pair of singularities of $f(x,z)$ that orbit the circle $|z| = \exp(|\Im(x_j)|)$ in opposite directions. The singularity at $x = 0$ similarly generates a pair of singularities on the unit circle at $z = \exp(\pm ix)$ in the complex $z$-plane.

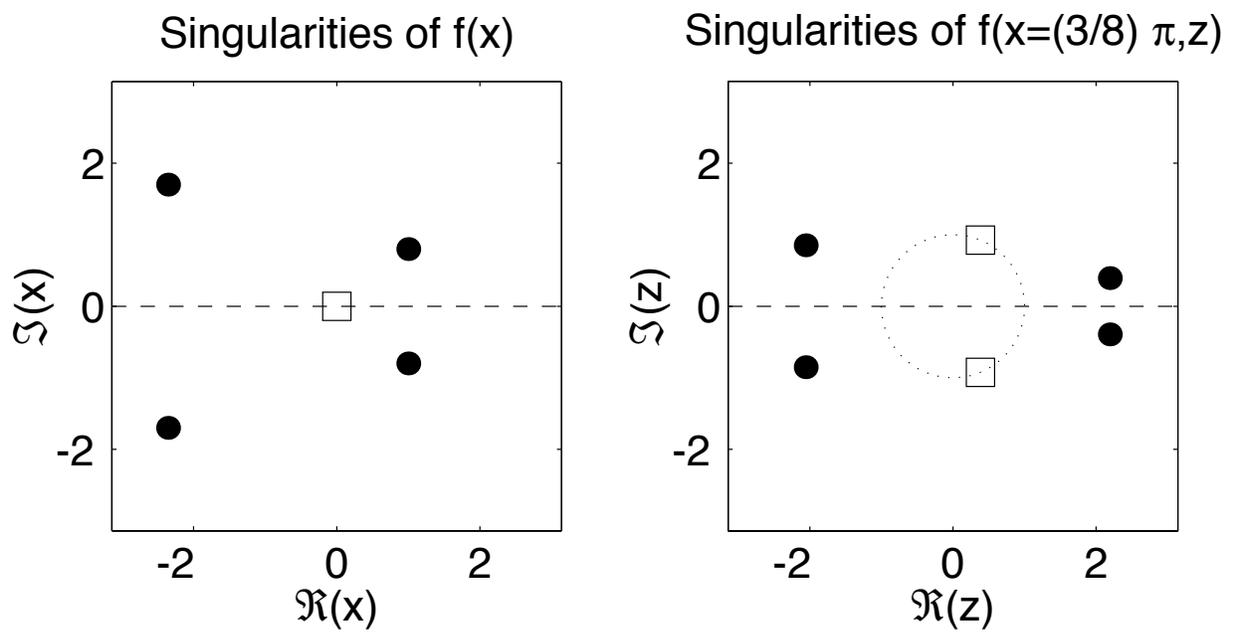

Figure 2: Schematic showing singularities of $f(x)$ on the left. For each real $x$, the Abel extension $f(x,z)$ has corresponding singularities in the $z$-plane; a typical case for the arbitrarily chosen $x = (3/8)\pi$ is shown on the right. The singularity on the real axis is marked with a square; the singularity off the real axis and its images are disks. The unit circle in the $z$-plane is shown as the dotted circle on the right.

Now that we know where the singularities of $f(x)$ translate into singularities of the Abel extension $f(x, z)$ in the $z$-plane when $x$, the argument of the bivariate extension, is real, we can address the question: where does the conformal map move these singularities in the $\zeta$-plane?

**Theorem 5** *If $f(z)$ has a singularity at $z = r_j \exp(i\theta_j)$, then $f(z(\zeta))$ aas a singularity at $\zeta_j$ where*

$$|\zeta_j| = \frac{2r_j}{\sqrt{1 + r_j^2 + 2r_j \cos(\theta_j)}} \tag{46}$$

*If the corresponding singularity of $f(x)$ is at $x = \sigma_j + i\tau_j$, then*

$$|\zeta_j(x)| = \frac{2\exp(|\tau_j|)}{\sqrt{1 + \exp(2|\tau_j|) + 2\exp(|\tau_j|) \cos(x - \sigma_j)}} \tag{47}$$

Proof: Let $\overline{z}$ denotes the complex conjugate of $z$. Recall that $\zeta = 2z/(1 + z)$. Then multiplying both sides by $\overline{\zeta}$ to obtain

$$|\zeta|^2 = \frac{4|z|^2}{1 + z + \overline{z} + |z|^2} \tag{48}$$

Substituting $z = r_j \exp(i\theta_j)$ and simplifying $z + \overline{z}$ to $2\cos(\arg(z))$ gives the theorem. QED

Thus, the magnitude of a singularity in $\zeta$ varies with the parameter $x$. However, some obvious remarks are possible. Because $|\zeta_j(x)|$ depends only on $x$ through its denominator, it follows that the smallest value of $|\zeta_j(x)|$ occurs when the denominator is a maximum, which happens when the argument of the cosine is zero:

$$\min_{x \in [-\pi, \pi]} |\zeta_j(x)| = \frac{2\exp(|\tau_j|)}{1 + \exp(|\tau_j|)} \tag{49}$$

**Theorem 6** *Let $f(x)$ be a $2\pi$-periodic analytic function which is singular only at $x = 0$ and $x_j = \sigma_j + i\tau_j$. Let $f_N^\sigma(x; N)$ denote its Euler-accelerated partial sum:*

$$f_N^\sigma \equiv \sum_{n=-N}^{N} \mu_n^N c_n \exp(inx) \tag{50}$$

*where $\mu_n^N$ denotes the weights of the Euler acceleration. Then*

$$|f(x) - f_N^\sigma(x)| \le p(x) \exp(-q(x)N) \tag{51}$$

*where*

$$\exp(q(x)) = \rho(x) = \begin{cases} \min(2, |\zeta_1(x)|, |\zeta_2(x)|, \dots), |x| \ge (2/3)\pi \\ \min(|\zeta_0(x)|, |\zeta_1(x)|, \dots), |x| < (2/3)\pi \end{cases} \tag{52}$$

*where*

$$|\zeta_j| = \frac{2r_j}{\sqrt{1 + r_j^2 + 2r_j \cos(\theta_j)}} \tag{53}$$

*and the image of the singularity at $x = 0$ is*

$$|\zeta_0(x)| = \frac{1}{\cos(x/2)} \tag{54}$$

*Furthermore, if the singularity on the real axis is no stronger than a Dirac delta-function, the error proportionality factor $p(x)$ may be bounded by a constant divided by $1/|x|$:*

$$|p(x)| \le C/|x| \tag{55}$$

Proof: The singularity at $\zeta = 2$ comes from the "metric" singularity in the map itself $z = \zeta/(2 - \zeta)$. The magnitude of the location in $\zeta$ of the other singularities is given by Theorem 5. The bound on the error in the Eulerized partial sum follows form applying the power series rate of convergence Theorem 1. Note the error estimatescomes entirely from the analysis of the power series in $\zeta$; this is applied to eah value of $x$ separately with $x$ reduced to merely a parameter.

To prove that the constant $C$ in the error bound (55) is independent of $x$, it is helpful to examine a particular case which bounds all likely cases of interest. The periodized Dirac delta function is

$$\delta(x) = \sum_{n=-\infty}^{\infty} \exp(inx) \tag{56}$$

Its Abel Extension is

$$\delta^{Abel}(x;z) = \sum_{n=-\infty}^{\infty} z^n \exp(inx) \tag{57}$$

$$= \lambda(x;z,0) \tag{58}$$

$$= \frac{(1-z^2)}{(1+z^2) - 2z\cos(x)} \tag{59}$$

$$= \frac{1}{1-\exp(ix)z} + \frac{1}{1-\exp(-ix)z} - 1 \tag{60}$$

where the last line follows by applying the geometric series to both fractions, and comparing the result with the first line. The conformal mapping $z = \zeta/(2 - \zeta)$ yields

$$\delta^{Abel}(x;\zeta) = \frac{1-\zeta/2}{1-(1+\exp(ix))\zeta/2} + \frac{1-\zeta/2}{1-(1+\exp(-ix))\zeta/2} - 1 \tag{61}$$

It is clear that the fractions may be again be expanded in geometric series:

$$\delta^{Abel}(x;\zeta) = (1-\zeta/2)\left\{\sum_{n=0}^{\infty} ((1+exp(ix))^n + (1+\exp(-ix))^n)\frac{\zeta^n}{2^n}\right\} - 1 \tag{62}$$

The error in truncating after the $N$-th term aat $\zeta = 1$ is

$$|\delta(x) - \delta_N^\sigma(x)| = \left(\frac{1+\exp(ix)}{2}\right)^{N+1} \frac{1}{1-(1+\exp(ix))/2}$$
$$+ \left(\frac{1+\exp(-ix)}{2}\right)^{N+1} \frac{1}{1-(1+\exp(-ix))/2} \tag{63}$$

It follows from this that proportionality constant $p(x)$ can be become unbounded as $|x| \to 0$, but no faster than $1/|x|$ as follows from Taylor expanding the denominators in (63). It is straightforward to extend this reasoning to other species of singularities at $x = 0$. QED

# 5 Interpretation

If $f(x)$, like the shifted sawtooth function, has singularities only at $x = 0$, or only has singularities at the origin and at locations far from the real $x$ axis, then only $\zeta_0(x)$, the image of the singularity at $x = 0$ and the metric singularity at $\zeta = 2$ are significant. Nothing that

$$|\zeta_0(x)| = \frac{\sqrt{2}}{\sqrt{1 + \cos(x)}} = \frac{1}{\cos(x/2)} \tag{64}$$

it follows that for the sawtooth function and others in the same class of no other singularities near the real axis,

$$\rho(x) = \min\left(2, \frac{1}{\cos(x/2)}\right) \tag{65}$$

$$\approx 1 + \frac{1}{8}x^2, \qquad |x| << 1 \tag{66}$$

which is graphed in Fig. 3. The terms in filtered Fourier series of degree $N$ converge like $(1/\rho(x)^N)$; when $|x| > \pi/2$, the metric sinugularity dominates and each additional term reduces the error by a factor of two. Unfortunately, the rate of convergence becomes slower and slower as we approach the location of the singularity of $f(x)$ on the real axis at $x = 0$.

In applications, most functions are more complicated than the sawtooth function in that they also have singularities off the real axis. A useful generalizatio is to add to the sawtooth function the "Symmetric Imbricated Lorentzian" [8]

$$\lambda(x; p, \phi) \equiv \frac{(1 - p^2)}{(1 + p^2) - 2p\cos(x - \phi)}$$

$$= 1 + 2\sum_{n=1}^{\infty} p^n \cos(n[x - \phi]) \tag{67}$$

$$= 2a \sum_{m=-\infty}^{\infty} \frac{1}{(\tau^2 + (x - 2\pi m - \phi)^2)} \tag{68}$$

where $\tau \equiv -\log(p) \to p = \exp(-\tau)$. The third form, which is known variously as the "periodization" or "imbrication" of the rational function, shows explicitly that $\lambda(x, p, \phi)$ is periodic and singularity-free on the real axis; when added to the shifted sawtooth function, the only singularity on the real interval $x \in [-\pi, \pi]$ is at $x = 0$. However, the imbricate form shows explicitly that $\lambda(x; p, \phi)$ has first order poles at $x = i \pm a + 2\pi m$ for any integer $m$. By varying the parameters $p$ and $\phi$, we can move these conjugate poles to any location that we please. Although the off-axis singularities in applications may be logarithms, square roots, higher order poles or something more exotic, this function with first order poles is much more representative (and therefore much more useful) than it might at first seem. The reason is that the convergence of spectral series depends, as discussed in Chapter 2 of [8], depends *exponetially* on the *location* of the singularities of the function defined by the series but only *algebraically* (i. e., very slowly) on the *type* of singularity. (The mantra of real estate, "Location! Location! Location!", applies to power series and Fourier series even more.)

The Abel extension of $\lambda$ is

$$\lambda(z(x); p, \phi) \equiv \frac{(1 - p^2 z^2)}{(1 + p^2 z^2) - 2pz\cos(x - \phi)} \tag{69}$$

The conformal mapping yields

$$\lambda(\zeta(x); p, \phi) \equiv \frac{(1 - p^2\zeta^2/(2-\zeta)^2)}{(1 + p^2\zeta^2/(2-\zeta)^2) - 2\,p[\zeta/(2-\zeta)]\,\cos(x-\phi)} \quad (70)$$

The Euler acceleration has long been known to be a "regular" acceleration in the sense that it can never stop a convergent series from converging, as shown by the analysis above. However, regular filters generally can *slow* the convergence of geometrically convergent series. If a function has an off-axis singularity at $\Im(x) = \tau$, then the coefficients of the unaccelerated series will be smaller than its predecessor by a factor of $\exp(|\tau|)$, which is equivalent to the radius of convergence of the Eulerized power series in $\zeta$. The rate of convergence of the Eulerized series varies with $x$ as already noted, but has a minimum at

$$\rho_{min} = \frac{2\exp(|\tau_j|)}{1 + \exp(|\tau_j|)} \approx 1 + |\tau/2|, \tau \ll 1 \quad (71)$$

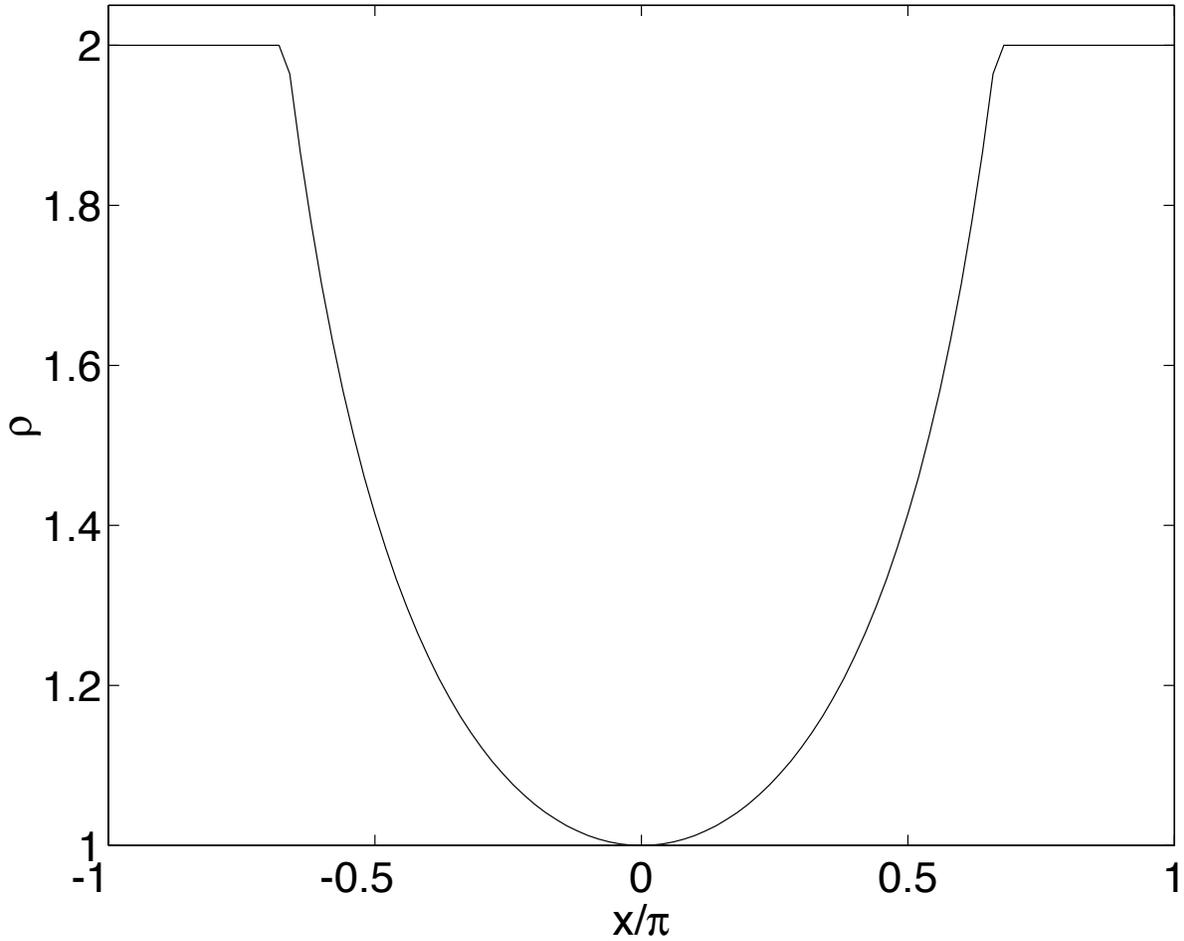

Figure 3: Plot of the radius of convergence $\rho(x)$ of the Euler-accelerated series for the sawtooth function for real $x$.

versus $\exp(|\tau|) \approx 1+|\tau|$ in the same limit. Thus, we need roughly twice as many Fourier terms (at the most unfavorable value of $x$) as for the unaccelerated series. Fig. 4 shows the convergence factor $\rho(x)$ for a typical case.

Off-axis singularities have the most pernicious effect on the rate of convergence when they are as far as possible (in $\Re(x)$) from the singularity on the real axis at $x = 0, \pm 2\pi, \pm 4\pi \ldots$. We therefore choose the phase factor $\phi$ so that there are poles on the line $\Re(x) = \pi$. As shown in Fig. 4, the Euler "acceleration" halves the convergence factor in the vicinity of $x = \pi$ if the off-axis singularities are only a short distance $\tau$ from the real axis.

Fig. 5 shows how the radius of convergence $\rho$ varies with and $x$ and with $\tau$, the distance of the poles from the real axis.

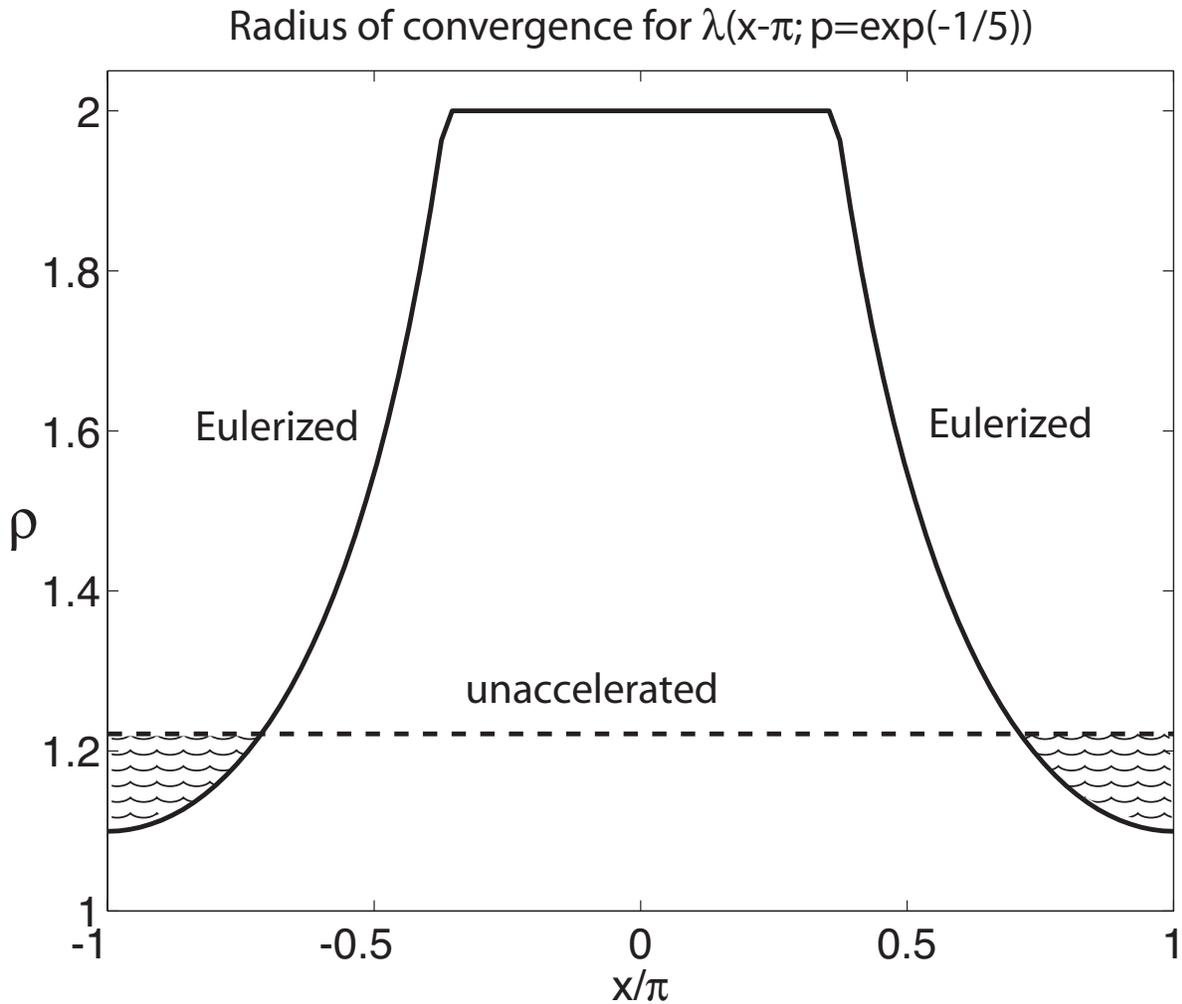

Figure 4: Plot of the radius of convergence $\rho(x)$ of the Euler-accelerated series for the function $\lambda(x - \pi, p = \exp(-1/5))$, which is singularity-free on the real axis, but has poles at $x = \pm i/5$. Shading marks the subintervals in $x$ where acceleration has slowed the radius of convergence.

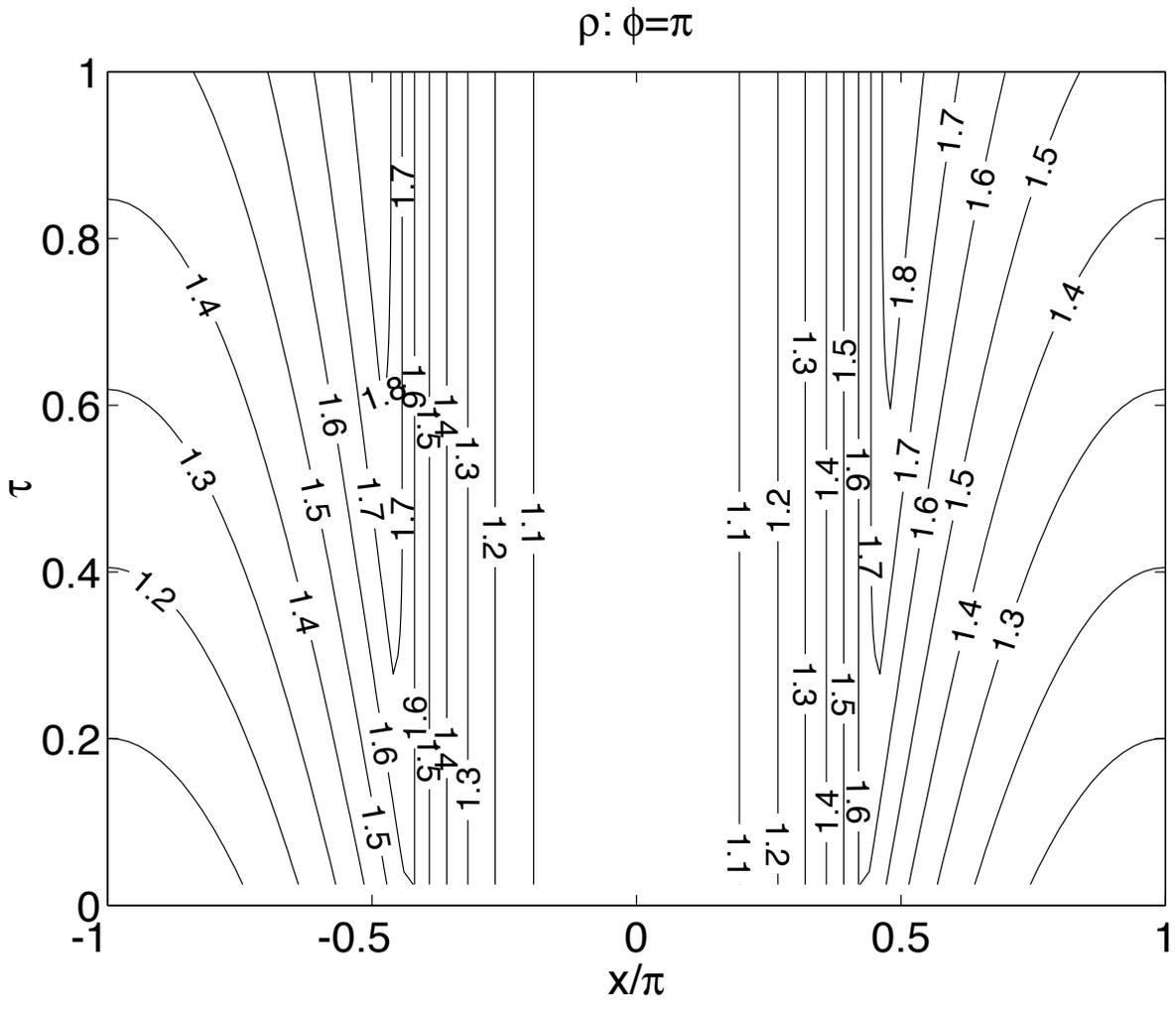

Figure 5: Plot of the radius of convergence $\rho(x)$ of the Euler-accelerated series for the function $f(x) = \text{Sws}(x) + \lambda(x - \pi, p = \exp(-\tau))$, which has poles at $x = \pm\tau$.

# 6 Numerical Illustration of Geometric Convergence

Fig. 7 shows the errors for the Euler-accelerated Fourier series. On a log-linear plot, a geometric rate of convergence appears as a straight line. The errors in the Fourier series fluctuate with $N$, but the *envelope* displays the predicted rate of convergence where the envelope is formally defined as follows.

**Definition 6 (Envelope of the Spectral Coefficients)**
*The ENVELOPE of the spectral coefficients is defined to be that curve, of the form of the leading asymptotic term in the* logarithm *of the absolute value of the spectral coefficients, which bounds the coefficients from above as tightly as possible.*

The graph features a dotted line which shows the predicted constant $\exp(-qN)/N$ where $q(x) = -\log(\cos(x/2))$; this is a very good approximation to the envelope.

Fig. 6 is identical except that $x$ is now only one-fifth as far from the jump discontinuity. The slope is much smaller, but the same geometric rate of convergence is evident.

# 7 Other Filters

We must point out that the Euler filter is not the optimum sum acceleration strategy, but only the filter that allows the most elementary proof. Fig. 8 compares the Euler filter (top) with two better filters: the Erfclog filter of Boyd [6] and the HDAF filter of Tanner [60]. The Erfclog filter replaces the Euler weights by

$$\sigma_{Erfc-Log}(\theta; p) \equiv \frac{1}{2}\mathrm{erfc}\left\{2p^{1/2}\overline{\theta}\sqrt{\frac{-\log(1-4\overline{\theta}^2)}{4\overline{\theta}^2}}\right\}, \qquad \overline{\theta} \equiv |\theta| - \frac{1}{2} \tag{72}$$

where for optimum results the "order" $p$ is

$$p(x) = 1 + N|x|/(2\pi) \tag{73}$$

Similarly, using Tanner's suggested optimum parameters of $\alpha = 1$ and $\kappa = 1/15$,

$$\sigma_{HDAF}(\theta) = \exp(-N|x|\theta^2/2) \sum_{j=0}^{\mathrm{floor}(N|x|/15)} (N|x|\theta^2/2)^j/j! \tag{74}$$

The guidelines show that all three display a geometric rate of convergence, even though $x$ is rather close to the jump discontinuity of the function $f(x)$. However, the Erfclog and HDAF filters are much superior to the Euler acceleration not only for the $x$ shown, but for all $x$. The HDAF filter is not as good as the Erfclog filter far from the singularity, but since the rate of convergence is fast at such $x$, this is not as important as the acceleration near the discontinuity.

# 8 Summary

Using only a rational conformal mapping ("Moebius transformation") and the usual theory for the convergence of a power series in the complex plane, we have rigourized our earlier proof [5] by more careful bookeeping of singularities of $f(x)$ off the real axis. We show that the Euler series acceleration, though ancient, is able to recover a geometric rate of convergence in the sense that the error in the series, truncated after degree $N$, is proportional to $\exp(-q(x)N)$ where $q(x) = -\log(\cos(x/2))$ with the singularity located at $x = 0$.

Unfortunately, the Euler acceleration and all local filters, even nonlinear, fail in the neighborhood of the singularity. For the Euler filter, $q(x) \approx x^2/8$ for small $x$. A more uniform strategy is to apply a nonlocal filter such as Gibbs reprojection, or to use a filter only on a domain excluding the neighborhood of the singularity, and treat the singularity region by a polynomial based strategy [14, 15].

Acknowledgments. This work was supported by NSF Grants OCE998636, OCE0451951 and ATM 0723440. A short version of this work has appeared in [12].

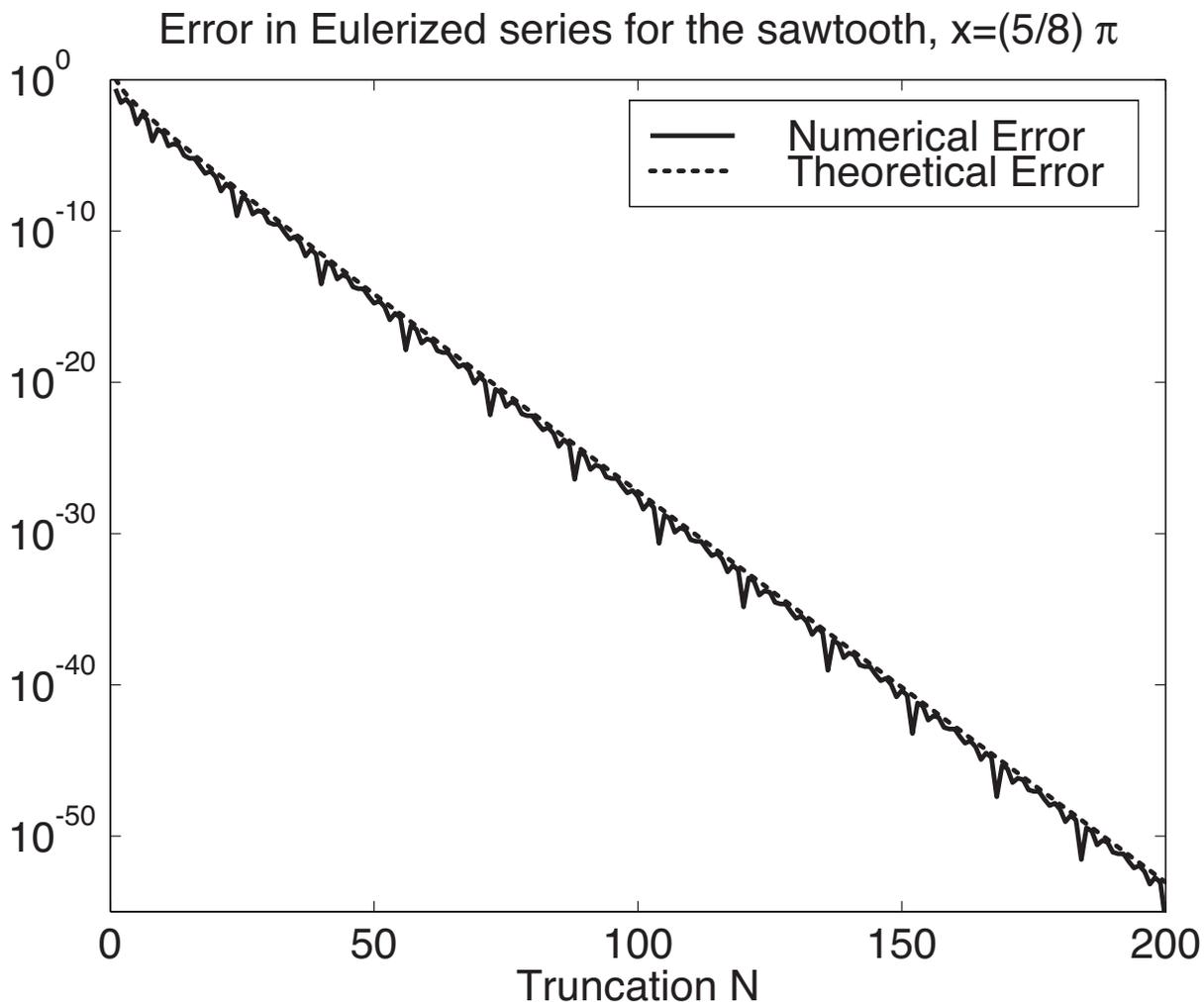

Figure 6: The solid curve is the error of the Euler-accelerated series for the function $f(x) = \text{Sws}(x)$, evaluated in Maple using 60 digit precision, for $x = (5/8)\pi$. The dotted line is the best fit to the envelope of the errors, $2\exp(-q(x)N)/N$ where $q(x) = -\log(\cos(x/2))$.

# References

[1] M. ABRAMOWITZ AND I. A. STEGUN, *Handbook of Mathematical Functions*, Dover, New York, 1965.

[2] R. ARCHIBALD, K. W. CHEN, A. GELB, AND R. RENAUT, *Improving tissue segmentation of human brain MRI through preprocessing by the Gegenbauer reconstruction method*, Neuroimage, 20 (2003), pp. 489-502.

[3] R. ARCHIBALD AND A. GELB, *A method to reduce the Gibbs ringing artifact in MRI scans while keeping tissue boundary integrity*, IEEE Trans. Medical Imaging, 21 (2002), pp. 305-319.

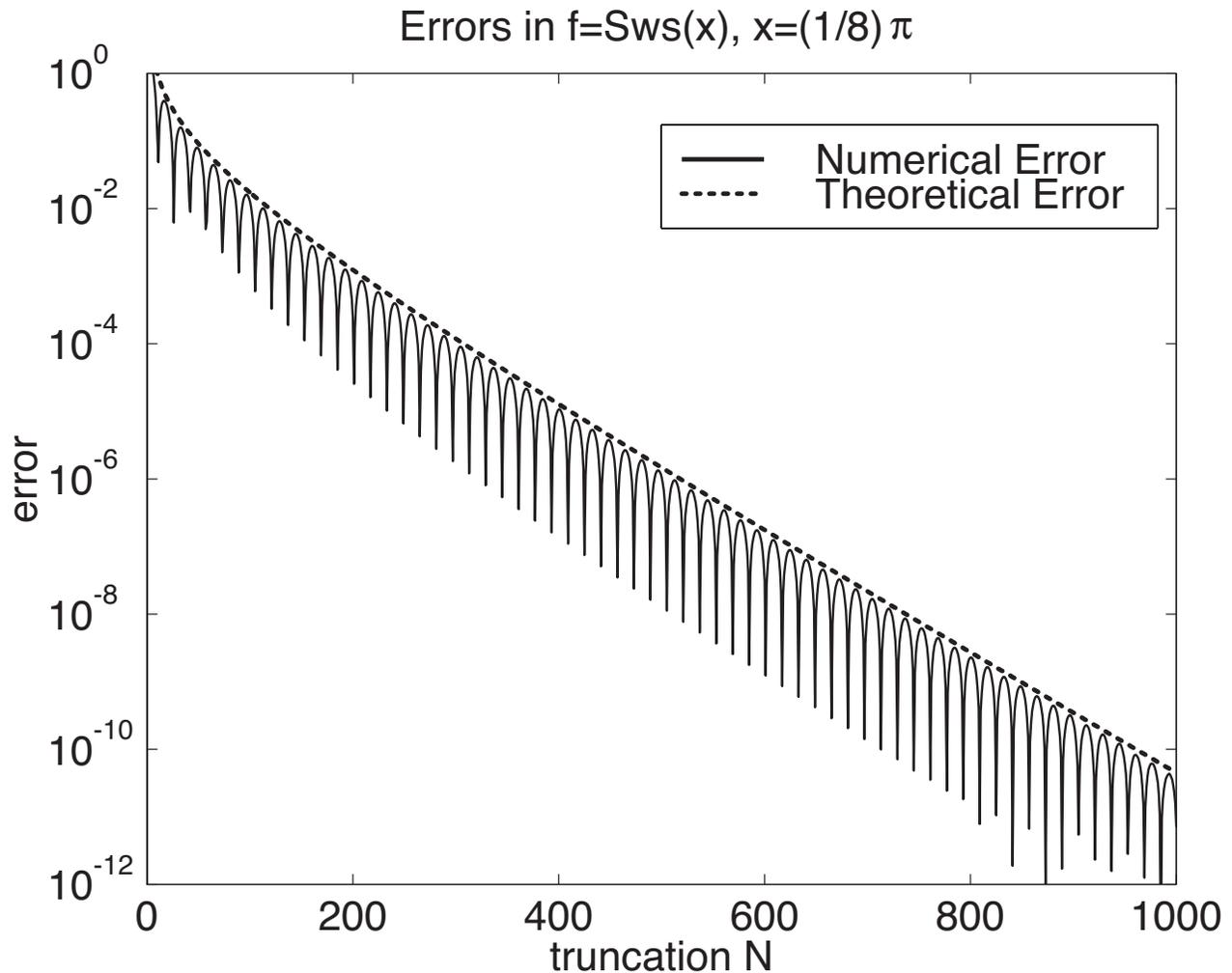

Figure 7: Same as previous figure but for $x = \pi/8$, much closer to the jump discontinuity of the shifted sawtooth function at $x = 0$. The dotted line is the best fit to the envelope of the errors, $12\exp(-q(x)N)/N$ where $q(x) = -\log(\cos(x/2))$.

[4] J. P. BOYD, *Sum-accelerated pseudospectral methods: The Euler-accelerated sinc algorithm*, App. Numer. Math., 7 (1991), pp. 287–296.

[5] ——, *A lag-averaged generalization of Euler's method for accelerating series*, Appl. Math. Comput., 72 (1995), pp. 146–166.

[6] ——, *The Erfc-Log filter and the asymptotics of the Vandeven and Euler sequence accelerations*, in Proceedings of the 3rd International Conference on Spectral and High Order Methods, A. V. Ilin and L. R. Scott, eds., Houston, 1996, Houston J. Mathematics, pp. 267–276.

[7] ——, *The Blasius function in the complex plane*, J. Experimental Math., 8 (1999), pp. 381–394.

[8] ——, *Chebyshev and Fourier Spectral Methods*, Dover, New York, 2001.

[9] ——, *Trouble with Gegenbauer reconstruction for defeating Gibbs' phenomenon: Runge phenomenon in the diagonal limit of Gegenbauer polynomial approximations*, J. Comput. Phys., 204 (2005), pp. 253–264.

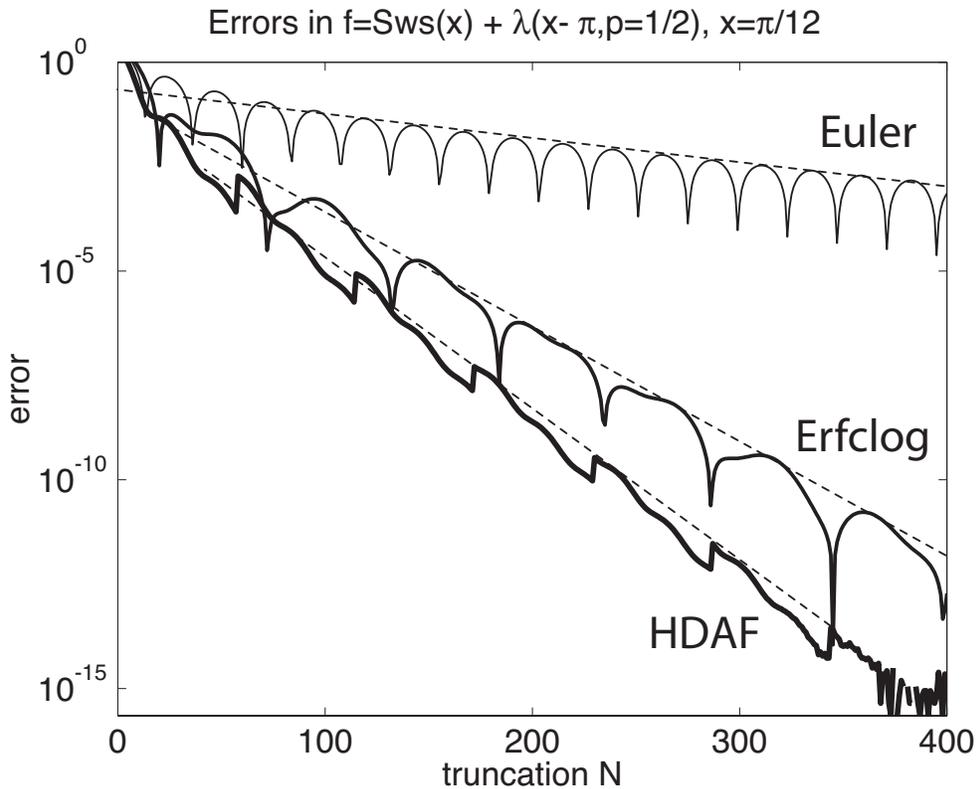

Figure 8: Errors $x = \pi/12$ versus the truncation $N$ of the Fourier series for $f(x) = \text{Sws}(x) + \lambda(x - \pi, p = 1/2)$ as accelerated by three different filters. The dashed guidelines show that the rate of convergence of the errors is proportional to $\exp(-qN)$ for some constant $q$ for all three filters. However, the constant $q$ is much better for the Erfclog and HDAF filters.

[10] ———, *Exponentially accurate Runge-Free approximation of non-periodic functions from samples on an evenly-spaced grid*, Appl. Math. Lett., 20 (2007), pp. 971–975.

[11] ———, *Acceleration of algebraically-converging Fourier series when the coefficients have series in powers of* $1/n$, J. Comput. Phys., 228 (2008), pp. 1401–1411.

[12] ———, *A proof, based on the euler sum acceleration, of the recovery of am exponential (geometric) rate of convergence for the Fourier series of a function with Gibbs Phenomenon*, in Proceedings of ICOSAHOM '09, E. Ronquist, ed., Lecture Notes in Computational Science and Engineering, New York, 2010, Springer-Verlag.

[13] J. P. BOYD AND D. W. MOORE, *Summability methods for Hermite functions*, Dyn. Atmos. Oceans, 10 (1986), pp. 51–62.

[14] J. P. BOYD AND J. R. ONG, *Exponentially-convergent strategies for defeating the Runge Phenomenon for the approximation of non-periodic functions, Part I: Single-Interval schemes*, Commun. Comput. Phys., 5 (2009), pp. 484–497.

[15] ———, *Exponentially-convergent strategies for defeating the runge phenomenon for the approximation of non-periodic functions, Part II: Multi-Interval and multi-approximation schemes*, Commun. Comput. Phys., (2009). to be submitted.

[16] O. P. BRUNO AND F. REITICH, *Approximation of analytic functions — a method of enhanced convergence*, Quart. J. Mech. Appl. Math., 63 (1994), pp. 195–213.

[17] W. CAI, D. GOTTLIEB, AND C. W. SHU, *On one-sided filters for spectral Fourier approximation of discontinuous functions*, SIAM J. Numer. Anal., 29 (1992), pp. 905–916.

[18] C. CANUTO, M. Y. HUSSAINI, A. QUARTERONI, AND T. A. ZANG, *Spectral Methods: Fundamentals in Single Domain*, Springer-Verlag, New York, 2006. 558 pp.

[19] T. A. DRISCOLL AND B. FORNBERG, *A Padé-based algorithm for overcoming the Gibbs phenomenon*, Numerical Algorithms, 26 (2001), pp. 77–92.

[20] K. S. ECKHOFF, *Approximate and efficient reconstruction of discontinuous functions from truncated series expansions*, Mathematics of Computation, 61 (1993), pp. 745–763.

[21] ———, *On discontinuous solutions of hyperbolic equations*, in Analysis, Algorithms and Applications of Spectral and High Order Methods for Partial Differential Equations, C. Bernardi and Y. Maday, eds., Selected Papers from the International Conference on Spectral and High Order Methods (ICOSAHOM '92), Le Corum, Montpellier, France, 22-26 June 1992, North-Holland, Amsterdam, 1994, pp. 103–112. Also in Comput. Meth. Appl. Mech. Engrg., 116.

[22] ———, *Accurate reconstructions of functions of finite regularity from truncated Fourier series expansions*, Mathematics of Computation, 64 (1995), pp. 671–690.

[23] ———, *On a high order numerical method for solving partial differential equations in complex geometries*, J. Sci. Comput., 12 (1997), pp. 119–138.

[24] ———, *On a high order numerical method for functions with singularities*, Math. Comput., 67 (1998), pp. 1063–1087.

[25] K. S. ECKHOFF AND J. H. ROLFSNES, *On nonsmooth solutions of linear hyperbolic systems*, Journal of Computational Physics, 125 (1996), pp. 1–15.

[26] P. A. Frost and E. Y. Harper, *An extended Padé procedure for constructing global approximations from asymptotic expansions: An explication with examples*, SIAM Rev., 18 (1976), pp. 62-91.

[27] J. Geer and N. S. Banerjee, *Exponentially accurate approximations to piece-wise smooth periodic functions*, J. Scient. Comput., 12 (1997), pp. 253-287.

[28] A. Gelb and E. Tadmor, *Spectral reconstruction of piecewise smooth functions from their discrete data*, Model. Math. Anal. Numer., 36 (2002), pp. 155-175.

[29] D. Gottlieb, B. Gustafsson, and P. Forssen, *On the direct Fourier method for computer tomography*, IEEE Trans. on Medical Imaging, 19 (2000), pp. 223-232.

[30] D. Gottlieb and C.-W. Shu, *Resolution properties of the Fourier method for discontinuous waves*, in Analysis, Algorithms and Applications of Spectral and High Order Methods for Partial Differential Equations, C. Bernardi and Y. Maday, eds., Selected Papers from the International Conference on Spectral and High Order Methods (ICOSAHOM '92), Le Corum, Montpellier, France, 22-26 June 1992, North-Holland, Amsterdam, 1994, pp. 27-38. Also in Comput. Meths. Appl. Mech. Engrg., vol. 116.

[31] D. Gottlieb and C.-W. Shu, *On the Gibbs phenomenon and its resolution*, SIAM Rev., 39 (1997), pp. 644-668.

[32] D. Gottlieb, C.-W. Shu, A. Solomonoff, and H. Vandeven, *On the Gibbs phenomenon I: recovering exponential accuracy from the Fourier partial sum of a non-periodic analytic function*, J. Comput. Appl. Math., 43 (1992), pp. 81-98.

[33] A. J. Guttmann, *Asymptotic analysis of power-series expansions*, in Phase Transitions and Critical Phenomena, C. Domb and J. L. Lebowitz, eds., no. 13 in Phase Transitions and Critical Phenomena, Academic Press, New York, 1989, pp. 1-234.

[34] J. S. Hesthaven, S. Gottlieb, and D. Gottlieb, *Spectral Methods for Time-Dependent Problems*, Cambridge University Press, Cambridge, 2007. 284 pp.

[35] G. H. Hoffman, *Extension of perturbation series by computer — viscous flow between infinite rotating disks*, J. Comput. Phys., 16 (1974), pp. 240-258.

[36] H. H. H. Homeier, *An asymptotically hierarchy-consistent, iterative sequence transformation for convergence acceleration of Fourier series*, Numerical Algorithms, 18 (1998), pp. 1-30.

[37] C. Hunter and S. M. Lee, *The analytic structure of oseen flow past a sphere as a function of reynolds number*, SIAM J. Appl. Math., 46 (1986), pp. 978-999.

[38] A. J. Jerri, ed., *Gibbs Phenomenon*, Sampling Publishing, Potsdam, New York, 2010.

[39] W. B. Jones and G. Hardy, *Accelerating the convergence of trigonometric approximations*, Math. Comp., 24 (1970), pp. 547-560.

[40] J. Jung and B. D. Shizgal, *Generalization of the inverse polynomial reconstruction method in the resolution of the Gibbs phenomenon*, J. Comput. Appl. Math., 172 (2004), pp. 131-151.

[41] ———, *Inverse polynomial reconstruction method in the resolution of the Gibbs phenomenon*, J. Sci. Comput., 25 (2005), pp. 367-399.

[42] J. Jung and B. D. Shizgal, *On the numerical convergence with the inverse polynomial reconstruction method for the resolution of the Gibbs phenomenon*, J. Comput. Phys., 224 (2007), pp. 477–488.

[43] I. M. Longman, *The summation of power series and Fourier series*, J. Comput. Appl. Math., 12 & 13 (1985), pp. 447–457.

[44] ———, *The summation of series*, Appl. Numer. Math., 2 (1986), pp. 135–141.

[45] J. N. Lyness, *The calculation of trigonometric Fourier coefficients*, J. Comput. Phys., 54 (1984), pp. 57–73.

[46] A. Majda, J. McDonough, and S. Osher, *The Fourier method for nonsmooth initial data*, Math. Comput., 32 (1978), pp. 1041–1081.

[47] S. L. Marshall, *Convergence acceleration of Fourier series by analytical and numerical application of Poisson's formula*, J. Phys. A, 31 (1998), pp. 2691–2704.

[48] ———, *On the analytical summation of Fourier series and its relation to the asymptotic behavior of Fourier transforms*, J. Phys. A, 31 (1998), pp. 9957–9973.

[49] D. W. Moore, *Planetary-gravity waves in an equatorial ocean*, PhD dissertation, University, Division of Engineering and Applied Sciences, 1968. 207 pp.

[50] P. M. Morse and H. Feshbach, *Methods of Theoretical Physics, (two volumes)*, McGraw-Hill, New York, 1953.

[51] A. Nersessian and A. Poghosyan, *On a rational linear approximation of Fourier series for smooth functions*, J. Sci. Comput., 26 (2006), pp. 111–125.

[52] C. Oleksy, *A convergence acceleration method of Fourier series*, Comput. Phys. Commun., 96 (1996), pp. 17–26.

[53] F. W. J. Olver, *Asymptotics and Special Functions*, Academic, New York, 1974.

[54] S. Paszkowski, *Convergence acceleration of orthogonal series*, Numer. Algorithms, 47 (2008), pp. 35–62.

[55] C. J. Pearce, *Transformation methods in the analysis of series for critical properties*, Adv. Phys., 27 (1978), pp. 89–148.

[56] R. E. Scraton, *The practical use of the Euler transform*, BIT, (1989), pp. 356–360.

[57] J. K. Shaw, L. W. Johnson, and R. D. Riess, *Accelerating the convergence of eigenfunction expansions*, Math. Comp., 30 (1976), pp. 469–477.

[58] E. Tadmor and J. Tanner, *Adaptive mollifiers for high resolution recovery of piecewise smooth data from its spectral information*, Foundations Comput. Math., 2 (2002), pp. 155–189.

[59] ———, *Adaptive filters for piecewise smooth spectral data*, IMA J. Numer. Anal., 25 (2005), pp. 635–647.

[60] J. Tanner, *Optimal filter and mollifier for piecewise smooth spectral data*, Math. Comput., 75 (2006), pp. 767–790.

[61] M. Tasche, *Accelerating convergence of Fourier series*, Math. Nachr., 90 (1979), pp. 123–134. In German.

[62] M. VAN DYKE, *Extension of Goldstein's series for the Oseen drag of a sphere*, J. Fluid Mech., 44 (1970), pp. 365–372.

[63] ———, *Computer extension of perturbation series in fluid mechanics.*, SIAM J. Appl. Math., 28 (1975), pp. 720–734.

[64] ———, *Perturbation Methods in Fluid Mechanics*, Parabolic Press, Stanford, California, 2d ed., 1975.

[65] ———, *Computer extended series*, Ann. Rev. Fluid Mech., 16 (1984), pp. 287–309.

[66] H. VANDEVEN, *Family of spectral filters for discontinuous problems*, J. Sci. Comput., 6 (1991), pp. 159–192.

[67] L. VOZOVOI, M. ISRAELI, AND A. AVERBUCH, *Analysis and application of the Fourier-Gegenbauer method to stiff differential equations*, SIAM J. Numer. Anal., 33 (1996), pp. 1844–1863.

[68] L. VOZOVOI, A. WEILL, AND M. ISRAELI, *Spectrally accurate solution of non-periodic differential equations by the Fourier-Gegenbauer method*, SIAM J. Numer. Anal., 34 (1997), pp. 1451–1471.

[69] S. WEINBERG, *Perturbation theory for strong repulsive potentials*, J. Math. Phys., 5 (1964), pp. 743–.

[70] R. K. WRIGHT, *A robust method for accurately representing nonperiodic functions given Fourier coefficient information*, J. Comput. Appl. Math., 140 (2002), pp. 837–848.

[71] ———, *Local spline approximation of discontinuous functions and location of discontinuities, given low-order Fourier coefficient information*, J. Comput. Appl. Math., 164–165 (2004), pp. 783–795.